\documentclass[twoside,11pt]{amsart}

\usepackage[cmtip,all]{xy}
\usepackage{xspace, enumerate, times, color}
\usepackage{amssymb, hyperref}

\usepackage{graphicx}
\input xypic
\input xy
\xyoption{all}
\def\sqr#1#2{{\vcenter{\hrule height.#2pt
        \hbox{\vrule width.#2pt height#1pt \kern#1pt
                \vrule width.#2pt}
        \hrule height.#2pt}}}

\numberwithin{equation}{section}

\newtheorem{theorem}{Theorem}[section]
\newtheorem{lemma}[theorem]{Lemma}

\newtheorem{open Problem}[theorem]{Open Problem}

\newtheorem{remark}[theorem]{Remark}

\newcommand{\be}{\begin{equation*}}
\newcommand{\ee}{\end{equation*}}
\newcommand{\bee}{\begin{equation}}
\newcommand{\eee}{\end{equation}}

\definecolor{lighterorange}{cmyk}{0,0.42,0.66,0.0}

\title[Partitions of finite nonnegative integer sets with identical representation functions]{Partitions of finite nonnegative integer sets with identical representation functions}
\author{Cui-Fang Sun,\;\; Hao Pan}

\begin{document}

\date{2022-8-12\\E-mail:  cuifangsun@163.com,\;\;  panhao0905@163.com}

\maketitle

\begin{abstract}
Let $\mathbb{N}$ be the set of all nonnegative integers. For $S\subseteq \mathbb{N}$ and $n\in \mathbb{N}$, let the representation function $R_{S}(n)$ denote the number of solutions of the equation $n=s+s'$ with $s, s'\in S$ and $s<s'$. In this paper, we determine the structure of $C, D\subseteq \mathbb{N}$ with $C\cup D=[0, m]$ and $|C\cap D|=2$ such that $R_{C}(n)=R_{D}(n)$ for any nonnegative integer $n$.

\noindent{{\bf Keywords:}\hspace{2mm} Representation function, Thue-Morse sequence, S\'{a}rk\"{o}zy's problem.}
\end{abstract}

\maketitle

\section{Introduction}
Let $\mathbb{N}$ be the set of all nonnegative integers. For $S\subseteq \mathbb{N}$ and $n\in \mathbb{N}$, let the representation function $R_{S}(n)$ denote the number of solutions of the equation $s+s'=n$ with $s, s'\in S$ and $s<s'$. Let $A$ be the set of all nonnegative integers which contain an even number of digits 1 in their binary representations and $B=\mathbb{N}\backslash A$. The set $A$ is called Thue-Morse sequence. For any positive integer $l$, let $A_{l}=A\cap [0, 2^{l}-1]$ and $B_{l}=B\cap [0, 2^{l}-1]$. For any nonnegative integer $r$ and $m$, let $r+m\mathbb{N}=\{r+mk: k\in\mathbb{N}\}$. S\'{a}rk\"{o}zy asked whether there exist two subsets $C, D\subseteq \mathbb{N}$ with $|(C\cup D)\backslash (C\cap D)|=\infty$ such that $R_{C}(n)=R_{D}(n)$ for all sufficiently large integers $n$. By using the Thue-Morse sequence, Dombi \cite{D} answered S\'{a}rk\"{o}zy's problem affirmatively.
In the last few decades, there are many results about representation functions (see, for example, \cite{CTY, CT, CB, KS1, JSYZ, L, LT, S, T, T1, TC, TL, YC}).

  In 2012, Yu and Tang \cite{YT} began to focus on partitions of nonnegative integers with the intersection not empty and obtain the following theorem.

\noindent {\bf Theorem A.} If $C\cup D=\mathbb{N}$ and $C\cap D=4\mathbb{N}$, then $R_{C}(n)=R_{D}(n)$ cannot hold for all sufficiently large $n$.

At the same time, they \cite{YT} posed the following conjecture.

\noindent {\bf Conjecture 1.} Let $m\in\mathbb{N}$ and $R\subset \{0, 1, \ldots, m-1\}$. If $C\cup D=\mathbb{N}$ and $C\cap D=\{r+km: k\in\mathbb{N}, r\in R\}$, then $R_{C}(n)=R_{D}(n)$ cannot hold for all sufficiently large $n$.

In 20l6, Chen and Lev \cite{CL} disproved the above conjecture by constructing a family of partitions of nonnegative integers.

\noindent {\bf Theorem B.} Let $l$ be a positive integer. There exist two sets $C$ and $D$ with $C\cup D=\mathbb{N}$ and $C\cap D=(2^{2l}-1)+(2^{2l+1}-1)\mathbb{N}$ such that $R_{C}(n)=R_{D}(n)$ for every positive integer $n$.

They also posed \cite{CL} the following two problems.

\noindent {\bf Problem 1.} Given $R_{C}(n)=R_{D}(n)$ for every positive integer $n$, $C\cup D=[0, m]$ and $C\cap D=\{r\}$ with $r\geq 0$ and $m\geq 2$, must there exist an integer $l\geq 1$ such that $r=2^{2l}-1, m=2^{2l+1}-2$, $C=A_{2l}\cup (2^{2l}-1+B_{2l})$ and $D=B_{2l}\cup (2^{2l}-1+A_{2l})$?

\noindent {\bf Problem 2.} Given $R_{C}(n)=R_{D}(n)$ for every positive integer $n$, $C\cup D=\mathbb{N}$ and $C\cap D=r+m\mathbb{N}$ with $r\geq 0$ and $m\geq 2$, must there exist an integer $l\geq 1$ such that $r=2^{2l}-1, m=2^{2l+1}-1$?

In 2017, Kiss and S\'{a}ndor \cite{KS} solved Problem 1 affirmatively. In 2021, Chen and Chen \cite{CC} solved Problem 2 affirmatively.

\noindent {\bf Theorem C.} Let $C$ and $D$ be sets of nonnegative integers such that $C\cup D=[0, m]$ and $|C\cap D|=1, 0\in C$. Then $R_{C}(n)=R_{D}(n)$ for every positive integer $n$ if and only if there exists a positive integer $l$ such that $C=A_{2l}\cup (2^{2l}-1+B_{2l})$ and $D=B_{2l}\cup (2^{2l}-1+A_{2l})$.

\noindent {\bf Theorem D.} Let $m\geq 2$ and $r\geq 0$ be two integers and let $C$ and $D$ be two sets with $C\cup D=\mathbb{N}$ and $C\cap D=r+m\mathbb{N}$ such that $R_{C}(n)=R_{D}(n)$ for every positive integer $n$. Then there exists a positive integer $l$ such that $r=2^{2l}-1$ and $m=2^{2l+1}-1$.

In this paper, we consider two sets $C$ and $D$ with $C\cup D=[0, m]$ and $|C\cap D|=2$ such that $R_{C}(n)=R_{D}(n)$ for any nonnegative integer $n$ and obtain the following results.

\begin{theorem}\label{thm1}
Let $m\geq 3$ be a positive integer. Let $C, D$ be two sets of nonnegative integers such that $C\cup D=[0, m]$, $C\cap D=\{r_{1}, r_{2}\}$, $r_{1}<r_{2}$ and $0\in C$.
If $r_{1}$ is even, then $R_{C}(n)=R_{D}(n)$ for any nonnegative integer $n$ if and only if there exists a positive integer $l$ such that $r_{1}=2^{2l+1}-2, r_{2}=2^{2l+1}-1, m=2^{2l+2}-3$ and
$$ C=A_{2l+1}\cup (2^{2l+1}-2+B_{2l+1}), \;\; D=B_{2l+1}\cup (2^{2l+1}-2+A_{2l+1}).$$
\end{theorem}

\begin{remark}
Let $m\geq 3$ be a positive integer. Let $C, D$ be two sets of nonnegative integers such that $C\cup D=[0, m]$, $C\cap D=\{r_{1}, r_{2}\}, r_{1}<r_{2}$ and $0\in C$. By Theorem $C$,
if $r_{1}$ is odd, then $R_{C}(n)=R_{D}(n)$ for any nonnegative integer $n$ if and only if there exists a positive integer $l$ such that $r_{1}=2^{2l}-1, r_{2}=2^{2l+1}+2^{2l}-2, m=2^{2l+2}-3$ and $$ C=A_{2l}\cup (2^{2l}-1+B_{2l})\cup \big(2^{2l+1}-1+(B_{2l}\cup (2^{2l}-1+A_{2l}))\big),$$
 $$ D=B_{2l}\cup (2^{2l}-1+A_{2l})\cup\big(2^{2l+1}-1+(A_{2l}\cup (2^{2l}-1+B_{2l}))\big).$$
\end{remark}

Throughout this paper, for $C, D\subseteq \mathbb{N}$ and $n\in \mathbb{N}$, let $R_{C, D}(n)$ be the number of solutions of $n=c+d$ with $c\in C$ and $d\in D$.
The characteristic function of $C\subseteq \mathbb{N}$ is denoted by
$$\chi_{C}(n)=\begin{cases}1,  &n\in C,\\
0, & n\not\in C.
\end{cases}$$

\section{Lemmas}

\begin{lemma}(\cite[Claim 1]{KS})\label{lem1}
Let $0< r_{1}< \cdots < r_{s}\leq m$ be integers. Then there exists at most one pair of sets $(C, D)$ such that $C\cup D=[0, m], 0\in C, C\cap D=\{r_{1}, \cdots, r_{s}\}$ and $R_{C}( n)=R_{D}(n)$ for every $n \leq m$.
\end{lemma}

\begin{lemma}(\cite[Claim 2]{KS})\label{lem2}
Let $(C, D)$ be a pair of different sets, $C\cup D=[0, m], C\cap D=\{r_{1}, \cdots, r_{s}\}$ and $R_{C}(n)=R_{D}(n)$ for every nonnegative integer $n$ and if $C'=m-C$ and $D'=m-D$, then $C'\cup D'=[0, m], C'\cap D'=\{m-r_{s}, \cdots, m-r_{1}\}$ and $R_{C'}(n)=R_{D'}(n)$ for every positive integer $n$.
\end{lemma}

\begin{lemma}(\cite[Claim 3]{KS})\label{lem3}
If for some positive integer M, the integers $M-1, M-2, M-4, M-8, \ldots, M-2^{\lceil \log_{2}M\rceil-1}$ are all contained in the set $A$, then $\lceil \log_{2}M\rceil$ is odd and  $M=2^{\lceil \log_{2}M\rceil}-1$.
\end{lemma}

\begin{lemma}(\cite[Claim 4]{KS})\label{lem6}
If for some positive integer M, the integers $M-1, M-2, M-4, M-8, \ldots, M-2^{\lceil \log_{2}M\rceil-1}$ are all contained in the set $B$, then $\lceil \log_{2}M\rceil$ is even and  $M=2^{\lceil \log_{2}M\rceil}-1$.
\end{lemma}

\begin{lemma}(\cite[Theorem 3]{KS})\label{lem4}
Let $C$ and $D$ be sets of nonnegative integers such that $C\cup D=[0, m], C\cap D=\emptyset$ and $0\in C$. Then $R_{C}(n)=R_{D}(n)$ for every positive integer $n$ if and only if there exists a positive integer $l$ such that $C=A_{l}$ and $D=B_{l}$.
\end{lemma}

\begin{lemma}(\cite[Corollary 1]{KS})\label{lem5}
If $C=A\cap [0, m]$ and $D=[0, m]$, where $m$ is a positive integer not of the form $2^{l}-1$, then there exists a positive integer $m<n<2m$ such that $R_{C}(n)\neq R_{D}(n)$.
\end{lemma}

\section{Proofs}

\begin{proof}[{\bf Proof of  Theorem 1.1}]
 We first prove the sufficiency of Theorem \ref{thm1}. Suppose that there exists a positive integer $l$ such that
 $$ C=A_{2l+1}\cup (2^{2l+1}-2+B_{2l+1}), \quad D=B_{2l+1}\cup (2^{2l+1}-2+A_{2l+1}).$$
 It is clear that $C\cup D=[0, 2^{2l+2}-3]$, $C\cap D=\{2^{2l+1}-2, 2^{2l+1}-1\}$ and $0\in C$. For any nonnegative integer $n$, we have
 \begin{eqnarray*}
 R_{C}(n)&=& R_{A_{2l+1}}(n)+R_{A_{2l+1}, 2^{2l+1}-2+B_{2l+1}}(n)+R_{2^{2l+1}-2+B_{2l+1}}(n)\\
  &=& R_{A_{2l+1}}(n)+ R_{A_{2l+1}, B_{2l+1}}(n-(2^{2l+1}-2))+R_{B_{2l+1}}(n-2(2^{2l+1}-2))
 \end{eqnarray*}
 and
 \begin{eqnarray*}
 R_{D}(n)&=& R_{B_{2l+1}}(n)+R_{2^{2l+1}-2+A_{2l+1}, B_{2l+1}}(n)+R_{2^{2l+1}-2+A_{2l+1}}(n)\\
  &=& R_{B_{2l+1}}(n)+ R_{A_{2l+1}, B_{2l+1}}(n-(2^{2l+1}-2))+R_{A_{2l+1}}(n-2(2^{2l+1}-2)).
 \end{eqnarray*}
By Lemma \ref{lem4}, for any nonnegative integer $k$, $R_{A_{2l+1}}(k)=R_{B_{2l+1}}(k)$. Thus $R_{C}(n)=R_{D}(n)$.

 Now we prove the necessity of Theorem \ref{thm1}. Suppose that $C\cup D=[0, m]$, $C\cap D=\{r_{1}, r_{2}\}, r_{1}<r_{2}, 2|r_{1}$ and $0\in C$ such that $R_{C}(n)=R_{D}(n)$ for any nonnegative integer $n$. Let $C'=m-C, D'=m-D$. By Lemma \ref{lem2}, we have $C'\cup D'=[0, m], C'\cap D'=\{m-r_{2}, m-r_{1}\}$ and $R_{C'}(n)=R_{D'}(n)$ for any positive integer $n$. Thus we may assume that $r_{1}+r_{2}\leq m$. An easy calculation shows that $r_{1}\geq 4$. Let
$$ p_{C}(x)=\sum\limits_{i=0}^{m}\chi_{C}(i)x^{i}, \;\; p_{D}(x)=\sum\limits_{i=0}^{m}\chi_{D}(i)x^{i}.$$
Then
\begin{equation}\label{3.2}
p_{D}(x)=\frac{1-x^{m+1}}{1-x}-p_{C}(x)+x^{r_{1}}+x^{r_{2}}.
\end{equation}
Since $R_{C}(n)=R_{D}(n)$ for any positive integer $n$, we have
$$\sum_{n=0}^{\infty} R_{C}(n)x^{n}=\sum_{n=0}^{\infty} R_{D}(n)x^{n}.$$
Then
\begin{equation}\label{3.3}
\frac{1}{2}(p_{C}(x)^{2}-p_{C}(x^{2}))=\sum_{n=0}^{\infty} R_{C}(n)x^{n}=\sum_{n=0}^{\infty} R_{D}(n)x^{n}=\frac{1}{2}(p_{D}(x)^{2}-p_{D}(x^{2})).
\end{equation}
By (\ref{3.2}) and (\ref{3.3}), we have
$$p_{C}(x)^{2}-p_{C}(x^{2})=\bigg(\frac{1-x^{m+1}}{1-x}-p_{C}(x)+x^{r_{1}}+x^{r_{2}}\bigg)^{2}-\bigg(\frac{1-x^{2m+2}}{1-x^{2}}-p_{C}(x^{2})+x^{2r_{1}}+x^{2r_{2}}\bigg).$$
Thus
\begin{eqnarray}\label{3.4}
2p_{C}(x^{2})&=& \frac{1-x^{2m+2}}{1-x^{2}}+2p_{C}(x)x^{r_{1}}+2p_{C}(x)x^{r_{2}}+2p_{C}(x)\frac{1-x^{m+1}}{1-x}\nonumber\\
&& -\bigg(\frac{1-x^{m+1}}{1-x}\bigg)^{2}-2x^{r_{1}}\frac{1-x^{m+1}}{1-x}-2x^{r_{2}}\frac{1-x^{m+1}}{1-x}-2x^{r_{1}+r_{2}}.
\end{eqnarray}

First, we will show that $r_{2}=2^{2l+1}-1$ and $r_{1}=2^{2l+1}-2$ for some positive integer $l$. If there exists a positive even integer $k$ such that $r_{1}\leq k < k+1< \min\{r_{2}, 2r_{1}\}\leq m$, then the coefficient of $x^{k}$ in (\ref{3.4}) is
\begin{equation}\label{3.5}
2\chi_{C}\bigg(\frac{k}{2}\bigg)=1+2\chi_{C}(k-r_{1})+2\sum\limits_{i=0}^{k}\chi_{C}(i)-(k+1)-2
\end{equation}
and the coefficient of $x^{k+1}$ in (\ref{3.4}) is
\begin{equation}\label{3.6}
0=2\chi_{C}(k+1-r_{1})+2\sum\limits_{i=0}^{k+1}\chi_{C}(i)-(k+2)-2.
\end{equation}
It follows from (\ref{3.5}) and (\ref{3.6}) that
\begin{equation}\label{3.7}
\chi_{C}\bigg(\frac{k}{2}\bigg)=\chi_{C}(k-r_{1})-\chi_{C}(k+1-r_{1})-\chi_{C}(k+1)+1.
\end{equation}
By Lemma \ref{lem1}, we have
\begin{equation}\label{3.8}
C\cap [0, r_{1}-1]=A\cap [0, r_{1}-1], \quad D\cap [0, r_{1}-1]=B\cap [0, r_{1}-1].
\end{equation}
Since $k+1-r_{1}<r_{1}, k-r_{1}$ is even, by (\ref{3.8}) and the definition of $A$, we have
$$\chi_{C}(k-r_{1})+\chi_{C}(k+1-r_{1})=1.$$

If $\chi_{C}(k-r_{1})=0$, then $\chi_{C}(k+1-r_{1})=1$. By (\ref{3.7}), we get $\chi_{C}(\frac{k}{2})=0$.

If $\chi_{C}(k-r_{1})=1$, then $\chi_{C}(k+1-r_{1})=0$. By (\ref{3.7}), we get $\chi_{C}(\frac{k}{2})=1$.

Thus $\chi_{C}(k-r_{1})=\chi_{C}\big(\frac{k}{2}\big)$. Putting $k=2r_{1}-2^{i+1}$ with $i\geq 0$, we have
$\chi_{C}(r_{1}-2^{i+1})=\chi_{C}(r_{1}-2^{i})$. Then
$$\chi_{C}(r_{1}-1)=\chi_{C}(r_{1}-2)=\chi_{C}(r_{1}-4)=\cdots=\chi_{C}(r_{1}-2^{\lceil\log_{2}r_{1}\rceil-1}).$$
By Lemmas \ref{lem3} and \ref{lem6}, we have $r_{1}=2^{\lceil\log_{2}r_{1}\rceil}-1$, which contradicts the condition that $r_{1}$ is even. Thus $r_{2}=r_{1}+1$. Let $k$ be a positive even integer with $r_{2}<k<k+1<2r_{1}<r_{1}+r_{2}\leq m$. Comparing the coefficients of $x^{k-1}, x^{k}$ and $x^{k+1}$ on the both sides of (\ref{3.4}) respectively, we have
$$ 0=2\chi_{C}(k-1-r_{1})+2\chi_{C}(k-1-r_{2})+2\sum\limits_{i=0}^{k-1}\chi_{C}(i)-k-4,$$
$$ 2\chi_{C}\bigg(\frac{k}{2}\bigg)=2\chi_{C}(k-r_{1})+2\chi_{C}(k-r_{2})+2\sum\limits_{i=0}^{k}\chi_{C}(i)-k-4,$$
$$ 0=2\chi_{C}(k+1-r_{1})+2\chi_{C}(k+1-r_{2})+2\sum\limits_{i=0}^{k+1}\chi_{C}(i)-k-6.$$
Calculating the above three equalities, we have
\begin{equation}\label{3.9}
\chi_{C}\bigg(\frac{k}{2}\bigg)=\chi_{C}(k-r_{1})-\chi_{C}(k-1-r_{2})+\chi_{C}(k).
\end{equation}
\begin{equation}\label{3.10}
\chi_{C}\bigg(\frac{k}{2}\bigg)=\chi_{C}(k-r_{2})-\chi_{C}(k+1-r_{1})-\chi_{C}(k+1)+1.
\end{equation}
By choosing $k=2r_{1}-2$ in (\ref{3.10}), we have
$$ 2\chi_{C}(r_{1}-1)=\chi_{C}(r_{1}-3)-\chi_{C}(2r_{1}-1)+1. $$
Then $\chi_{C}(r_{1}-1)=\chi_{C}(r_{1}-3)$. Thus $r_{1}\equiv 2\pmod 4$ and $r_{2}\equiv 3\pmod 4$.

If $k-1-r_{2}\equiv 0\pmod 4$, then $k-r_{1}\equiv 2\pmod 4$ and $k\equiv 0\pmod 4$. Thus
$$\chi_{C}\bigg(\frac{k-1-r_{2}}{2}\bigg)+\chi_{C}\bigg(\frac{k-r_{1}}{2}\bigg)=1.$$
Hence
$$\chi_{C}(k-1-r_{2})+\chi_{C}(k-r_{1})=1.$$
If $\chi_{C}(k-1-r_{2})=0$, then $\chi_{C}(k-r_{1})=1$. By (\ref{3.9}), we have $\chi_{C}(\frac{k}{2})=1$. If $\chi_{C}(k-1-r_{2})=1$, then $\chi_{C}(k-r_{1})=0$. By (\ref{3.9}), we have $\chi_{C}(\frac{k}{2})=0$. Thus $\chi_{C}\big(\frac{k}{2}\big)=\chi_{C}(k-r_{1})$ and $\chi_{C}\big(\frac{k}{2}\big)+\chi_{C}(k-1-r_{2})=1$. Noting that $\chi_{C}(k-1-r_{2})+\chi_{C}(k-r_{2})=1$, we have $\chi_{C}\big(\frac{k}{2}\big)=\chi_{C}(k-r_{2})$.

If $k-1-r_{2}\equiv 2\pmod 4$, then $k-r_{1}\equiv 0\pmod 4$ and $k\equiv 2\pmod 4$. By (\ref{3.9}), we have
$$ \chi_{C}\bigg(\frac{k-2}{2}\bigg)=\chi_{C}(k-2-r_{1})-\chi_{C}(k-3-r_{2})+\chi_{C}(k-2).$$
Then $\chi_{C}\big(\frac{k-2}{2}\big)=\chi_{C}(k-2-r_{1})$. Noting that $\chi_{C}\big(\frac{k-2}{2}\big)+\chi_{C}\big(\frac{k}{2}\big)=1$ and $\chi_{C}(k-1-r_{2})+\chi_{C}(k-r_{2})=1$, we have $\chi_{C}\big(\frac{k}{2}\big)=\chi_{C}(k-r_{2})$.

Thus we get $\chi_{C}\big(\frac{k}{2}\big)=\chi_{C}(k-r_{2})$. Let $k=2r_{2}-2^{i+1}$ with $i\geq 0$. Then
$\chi_{C}(r_{2}-2^{i+1})=\chi_{C}(r_{2}-2^{i})$. Thus
$$\chi_{C}(r_{2}-1)=\chi_{C}(r_{2}-2)=\chi_{C}(r_{2}-4)=\cdots=\chi_{C}(r_{2}-2^{\lceil\log_{2}r_{2}\rceil-1}).$$
By Lemmas \ref{lem3} and \ref{lem6}, we have $r_{2}=2^{\lceil\log_{2}r_{2}\rceil}-1$.

Comparing the coefficients of $x^{r_{1}}$ and $x^{r_{1}-1}$ on the both sides of (\ref{3.4}) respectively, we have
$$ 2\chi_{C}\bigg(\frac{r_{1}}{2}\bigg)=2\sum\limits_{i=0}^{r_{1}}\chi_{C}(i)-r_{1},\;\; 0=2\sum\limits_{i=0}^{r_{1}-1}\chi_{C}(i)-r_{1}.$$
Subtracting the above two equalities and dividing by 2, we get $\chi_{C}(\frac{r_{1}}{2})=1$.
Then $\chi_{A}(\frac{r_{1}}{2})=\chi_{C}(\frac{r_{1}}{2})=1$. Thus $r_{2}=2^{2l+1}-1$ and $r_{1}=2^{2l+1}-2$ for some positive integer $l$.

As $r_{1}+r_{2}=2^{2l+2}-3\leq m$, Lemma \ref{lem1} and the sufficiency of Theorem \ref{thm1}, we have
$$ C\cap [0, 2^{2l+2}-3]=A_{2l+1}\cup (2^{2l+1}-2+B_{2l+1}),\;\; D\cap [0, 2^{2l+2}-3]=B_{2l+1}\cup (2^{2l+1}-2+A_{2l+1}).$$

Next we will prove that $m=2^{2l+2}-3$. Now we prove that $m<3\cdot 2^{2l+1}-5$. Assume that $m\geq 3\cdot 2^{2l+1}-4$. Define the sets $E$ and $F$ by
\begin{eqnarray*}
&& E=A_{2l+1}\cup (2^{2l+1}-2+B_{2l+1})\cup (2^{2l+2}-2+(B_{2l+1} \cap [0, 2^{2l+1}-3])),\\
&& F=B_{2l+1}\cup (2^{2l+1}-2+A_{2l+1})\cup (2^{2l+2}-2+(A_{2l+1}\cap [0, 2^{2l+1}-3])).
\end{eqnarray*}
We will prove that $C\cap [0, 3\cdot 2^{2l+1}-5]=E$ and $D\cap [0, 3\cdot 2^{2l+1}-5]=F$.

If $n\leq 2^{2l+2}-3$, then $R_{E}(n)=R_{F}(n)$. If $2^{2l+2}-3<n\leq 3\cdot 2^{2l+1}-5$, then
\begin{eqnarray*}
R_{E}(n)&=& R_{A_{2l+1}, 2^{2l+1}-2+B_{2l+1}}(n)+R_{2^{2l+1}-2+B_{2l+1}}(n)\\
&& +R_{A_{2l+1}, 2^{2l+2}-2+(B_{2l+1}\cap [0, 2^{2l+1}-3])}(n)\\
&=& R_{A_{2l+1}, B_{2l+1}}(n-(2^{2l+1}-2))+R_{B_{2l+1}}(n-2(2^{2l+1}-2))\\
&&+R_{A_{2l+1}, B_{2l+1}}(n-(2^{2l+2}-2))
\end{eqnarray*}
and
\begin{eqnarray*}
R_{F}(n)&=& R_{B_{2l+1}, 2^{2l+1}-2+A_{2l+1}}(n)+R_{2^{2l+1}-2+A_{2l+1}}(n)\\
&& +R_{B_{2l+1}, 2^{2l+2}-2+(A_{2l+1}\cap [0, 2^{2l+1}-3])}(n)\\
&=& R_{B_{2l+1}, A_{2l+1}}(n-(2^{2l+1}-2))+R_{A_{2l+1}}(n-2(2^{2l+1}-2))\\
&&+R_{B_{2l+1}, A_{2l+1}}(n-(2^{2l+2}-2)).
\end{eqnarray*}
By Lemma \ref{lem4}, $R_{A_{2l+1}}(k)=R_{B_{2l+1}}(k)$ holds for any positive integers $k$ and then $R_{E}(n)=R_{F}(n)$. By Lemma \ref{lem1}, we have
$C\cap [0, 3\cdot 2^{2l+1}-3]=E$ and $D\cap [0, 3\cdot 2^{2l+1}-3]=F$.

We will prove that $3\cdot 2^{2l+1}-4\in C$ and $3\cdot 2^{2l+1}-3\in D$. By $3\cdot 2^{2l+1}-4=(2^{2l+1}-2)+(2^{2l+2}-2)$ in $D$, we have
$$ R_{C}(3\cdot 2^{2l+1}-4)=\chi_{C}(3\cdot 2^{2l+1}-4)+R_{B_{2l+1}}(2^{2l+1})+R_{A_{2l+1}, B_{2l+1}}(2^{2l+1}-2)$$
and
$$ R_{D}(3\cdot 2^{2l+1}-4)=1+R_{A_{2l+1}}(2^{2l+1})+R_{B_{2l+1}, A_{2l+1}}(2^{2l+1}-2).$$
Thus $R_{C}(3\cdot 2^{2l+1}-4)=R_{D}(3\cdot 2^{2l+1}-4)$ if and only if $\chi_{C}(3\cdot 2^{2l+1}-4)=1$, that is, $3\cdot 2^{2l+1}-4 \in C$. By
$$ R_{C}(3\cdot 2^{2l+1}-3)=\chi_{C}(3\cdot 2^{2l+1}-3)+R_{A_{2l+1}, B_{2l+1}}(2^{2l+1}-1)+R_{B_{2l+1}}(2^{2l+1}+1)$$
and
$$ R_{D}(3\cdot 2^{2l+1}-3)=R_{B_{2l+1}, A_{2l+1}}(2^{2l+1}-1)+R_{A_{2l+1}}(2^{2l+1}+1).$$
Thus $R_{C}(3\cdot 2^{2l+1}-3)=R_{D}(3\cdot 2^{2l+1}-3)$ if and only if $\chi_{C}(3\cdot 2^{2l+1}-3)=0$, that is, $3\cdot 2^{2l+1}-3 \in D$.

On the other hand, noting that $2^{2l+1}-2\in A_{2l+1}, 2^{2l+1}-1\in B_{2l+1}$, $3\cdot 2^{2l+1}-2=(2^{2l+1}-1)+(2^{2l+2}-1)$ in $C$ and $3\cdot 2^{2l+1}-2=1+(3\cdot 2^{2l+1}-3)$ in $D$, we obtain
\begin{eqnarray*}
 R_{C}(3\cdot 2^{2l+1}-2)&=&1+\chi_{C}(3\cdot 2^{2l+1}-2)+R_{A_{2l+1}, B_{2l+1}}(2^{2l+1})+R_{B_{2l+1}}(2^{2l+1}+2)\\
 && -\chi_{A_{2l+1}}(3\cdot 2^{2l+1}-2-(2^{2l+2}-2+2^{2l+1}-1))\\
 &=& 1+\chi_{C}(3\cdot 2^{2l+1}-2)+R_{A_{2l+1}, B_{2l+1}}(2^{2l+1})+R_{B_{2l+1}}(2^{2l+1}+2)
 \end{eqnarray*}
and
\begin{eqnarray*}
R_{D}(3\cdot 2^{2l+1}-2)&=& 1+R_{B_{2l+1}, A_{2l+1}}(2^{2l+1})+R_{A_{2l+1}}(2^{2l+1}+2)\\
 && -\chi_{B_{2l+1}}(3\cdot 2^{2l+1}-2-(2^{2l+2}-2+2^{2l+1}-2))\\
&=& R_{B_{2l+1}, A_{2l+1}}(2^{2l+1})+R_{A_{2l+1}}(2^{2l+1}+2).
\end{eqnarray*}
Thus by Lemma \ref{lem4}, we have $R_{C}(3\cdot 2^{2l+1}-2)>R_{D}(3\cdot 2^{2l+1}-2)$, which is impossible. Therefore $m< 3\cdot 2^{2l+1}-5$.
It follows that
\begin{eqnarray*}
&& C=A_{2l+1}\cup (2^{2l+1}-2+B_{2l+1})\cup (2^{2l+2}-2+(B_{2l+1} \cap [0, m-(2^{2l+2}-2)])),\\
&& D=B_{2l+1}\cup (2^{2l+1}-2+A_{2l+1})\cup (2^{2l+2}-2+A_{2l+1}\cap [0, m-(2^{2l+2}-2)])).
\end{eqnarray*}

We will show that $m=2^{2l+2}-3$. Assume that $m>2^{2l+2}-3$. If $m-(2^{2l+2}-2)\neq 2^{k}-1$ for any nonnegative integer $k$, then by Lemma \ref{lem5}, there exists an integer $u$ satisfying
$$m-(2^{2l+2}-2)<u<2(m-(2^{2l+2}-2))$$
such that
$$R_{A\cap [0, m-(2^{2l+2}-2)]}(u)\neq R_{B\cap [0, m-(2^{2l+2}-2)]}(u).$$
Noting that
$$m+(2^{2l+2}-2)<u+2(2^{2l+2}-2)<2m,$$
we have
\begin{eqnarray*}
R_{C}(2(2^{2l+2}-2)+u)&=& R_{2^{2l+2}-2+(B_{2l+1} \cap [0, m-(2^{2l+2}-2)])}(2(2^{2l+2}-2)+u)\\
&=&  R_{B_{2l+1} \cap [0, m-(2^{2l+2}-2)]}(u)
\end{eqnarray*}
and
\begin{eqnarray*}
R_{D}(2(2^{2l+2}-2)+u)&=& R_{2^{2l+2}-2+(A_{2l+1} \cap [0, m-(2^{2l+2}-2)])}(2(2^{2l+2}-2)+u)\\
&=&  R_{A_{2l+1} \cap [0, m-(2^{2l+2}-2)]}(u).
\end{eqnarray*}
By $R_{C}(2(2^{2l+2}-2)+u)=R_{D}(2(2^{2l+2}-2)+u)$, we have $R_{A_{2l+1} \cap [0, 2^{2l+2}-2]}(u)=R_{B_{2l+1} \cap [0, 2^{2l+2}-2]}(u)$, which is false. Thus we may assume that $m-(2^{2l+2}-2)=2^{k}-1$, where $k<2l+1$. Hence
\begin{eqnarray*}
&& C=A_{2l+1}\cup (2^{2l+1}-2+B_{2l+1})\cup (2^{2l+2}-2+B_{k}),\\
&& D=B_{2l+1}\cup (2^{2l+1}-2+A_{2l+1})\cup (2^{2l+2}-2+A_{k}).
\end{eqnarray*}
If $k=0$, then $C, D\subseteq [0, 2^{2l+2}-2]$ and $2^{2l+2}-4, 2^{2l+2}-2\in D$. Thus $R_{C}(2^{2l+3}-6)=0$ and $R_{D}(2^{2l+3}-6)=1$, a contradiction. Therefore we assume that $k>0$. Let $C'=2^{2l+2}+2^{k}-3-C$ and $D'=2^{2l+2}+2^{k}-3-D$. Now we need to consider the following two cases according to the parity of $k$.

{\bf Case 1.} $2\nmid k$. Then
\begin{eqnarray*}
&& C'=A_{k}\cup (2^{k}+A_{2l+1})\cup (2^{2l+1}+2^{k}-2+B_{2l+1}),\\
&& D'=B_{k}\cup (2^{k}+B_{2l+1})\cup (2^{2l+1}+2^{k}-2+A_{2l+1}).
\end{eqnarray*}

{\bf Case 2.} $2\mid k$. Then
\begin{eqnarray*}
&& C'=B_{k}\cup (2^{k}+A_{2l+1})\cup (2^{2l+1}+2^{k}-2+B_{2l+1}),\\
&& D'=A_{k}\cup (2^{k}+B_{2l+1})\cup (2^{2l+1}+2^{k}-2+A_{2l+1}).
\end{eqnarray*}

Thus
$$C'\cup D'=[0, 2^{2l+2}+2^{k}-3], \quad C'\cap D'=\{2^{2l+1}+2^{k}-2, 2^{2l+1}+2^{k}-1\}.$$
By Lemma \ref{lem2}, we know $R_{C'}(n)=R_{D'}(n)$ for all positive integers $n$. By Lemma \ref{lem1}, we need to examine the following two cases.

{\bf Case I.} $C'\cap [0, 2^{2l+1}+2^{k}-3]=A\cap  [0, 2^{2l+1}+2^{k}-3]$ and $D'\cap [0, 2^{2l+1}+2^{k}-3]=B\cap  [0, 2^{2l+1}+2^{k}-3]$.

{\bf Case II.} $C'\cap [0, 2^{2l+1}+2^{k}-3]=B\cap  [0, 2^{2l+1}+2^{k}-3]$ and $D'\cap [0, 2^{2l+1}+2^{k}-3]=A\cap  [0, 2^{2l+1}+2^{k}-3]$

If Case 1 and Case I hold, then $C'\cap [0, 2^{k}-1]=A\cap [0, 2^{k}-1]=A_{k}$. By the definition of $A$, we have $C'\cap [2^{k}, 2^{k+1}-1]=2^{k}+B_{k}$. However, by Case 1, we have $C'\cap [2^{k}, 2^{k+1}-1]=2^{k}+A_{k}$, a contradiction.

If Case 1 and Case II hold,  then $C'\cap [0, 2^{k}-1]=B\cap [0, 2^{k}-1]=B_{k}$. However, by Case 1, we have $C'\cap [0, 2^{k}-1]=A_{k}$, a contradiction.

If Case 2 and Case I hold, then $C'\cap [0, 2^{k}-1]=A\cap [0, 2^{k}-1]=A_{k}$. However, by Case 2, we have $C'\cap [0, 2^{k}-1]=B_{k}$, a contradiction.

If Case 2 and Case II hold, then by $2^{k}+1< 2^{2l+1}-1$ and the definition of $A$, we have $2^{k}+1\in A_{2l+1}$. Then $2^{k}+(2^{k}+1)\in 2^{k}+A_{2l+1}\subseteq C'$. By Case II, we have $2^{k}+(2^{k}+1)\in B$. However, $2^{k}+(2^{k}+1)=2^{k+1}+1\in A$, also a contradiction.

Therefore $m=2^{2l+2}-3$.

This completes the proof of Theorem \ref{thm1}.
\end{proof}

\end{document}